\documentclass[10pt,a4paper]{article}

\usepackage{latexsym}
\usepackage{amsfonts}
\usepackage{amssymb}
\usepackage{amsmath}
\usepackage{theorem}
\usepackage{enumerate}
\usepackage{cite}
\usepackage{url}
\usepackage[dvips]{color}
\usepackage{graphicx}

\theoremstyle{plain}
\theorembodyfont{\itshape}
\newtheorem{thm}{Theorem}
\newtheorem{lem}{Lemma}

\newtheorem{cor}{Corollary}

\theorembodyfont{\upshape}
\newtheorem{exem}{Example}

\newcommand{\proof}{\noindent {\bf Proof:} \hspace{0.1in}}
\newcommand{\qed}{\hfill\mbox{\raggedright $\Box$}\medskip}
\newcommand{\mydate}{
 \ifcase\month \or
 January\or February\or March\or April\or May\or June\or
 July\or August\or September\or October\or November\or December\fi
 \space \number\year}
\newcommand{\smin}{\,\raisebox{0.06em}{${\scriptstyle \in}$}\,}

\newcommand{\smwedge}{{\scriptstyle \wedge\,}}

\newcommand{\bwedge}{\raisebox{0.2ex}{${\textstyle \bigwedge}$}}

\newcommand{\bp}{\left(}
\newcommand{\ep}{\right)}
\newcommand{\ic}{\mathrm{i}}
\newcommand{\R}{\ensuremath{\mathbb{R}}}
\newcommand{\bdm}{\begin{displaymath}}
\newcommand{\edm}{\end{displaymath}}

\newcommand{\edem}{\hfill $\Box$ \\[3mm]}

\begin{document}

\title{Some maximal isotropic distributions and their relation to field theory }
\author{Leandro G. Gomes} 
\date{Instituto de Ci\^encias Exatas \\ Universidade Federal de Itajub\'a\\
      Itajub\'a, MG, Brazil \\ (lggomes@unifei.edu.br)}
      
\maketitle

\begin{abstract}
\noindent
 We study the behavior of differential forms in a manifold having at least one of
 their maximal isotropic local distributions endowed with the special algebraic 
 property of being decomposable. We show that they can be represented
 as the sum of a form with constant coefficients and one that vanishes whenever contracted
 with vector fields in the former distribution, provided some simple integrability 
 conditions are ensured. We also classify possible 'canonical coordinates' for
 a certain class of forms with potential applications in classical field theory.  
\end{abstract}

\section*{Introduction}

  The importance of the study of maximal isotropic (local) distributions with respect 
to a $(n+1)$-differential form $\omega$ on a manifold $P$ relies, at least but not last, 
in its connection with a covariant and finite dimensional approach to the classical 
theory of fields started with de Donder (\cite{DD}) and Weyl (\cite{We}). 
The existence of a special kind of such a 
distribution ensures that there are canonical coordinates 
$(x^\mu, q_i, p, p_i^\mu)$ in which $\omega$ locally emulates   
a preexisting canonical form $\Omega_0$ (\cite{tese, FG})), 
described by the formula
\begin{equation} \label{eq:MSPLF1}
  \Omega_0=~dq_{}^i \>\smwedge\, dp\>\!_i^\mu \>\smwedge\, d^{\,n} x_\mu^{} \, - \,
    dp \;\smwedge\, d^{\,n} x~,
\end{equation} 
where $d^{\,n} x=dx^{\,1} \smwedge \ldots \smwedge dx^{\,n}$ and 
$d^{\,n} x_\mu^{}~=~\mathrm{i}_{\partial_\mu^{}}^{} d^{\,n} x$.  
The local description of $\Omega_0$ given by equation (\ref{eq:MSPLF1})
permits the association of this geometric object to the de Donder-Weyl
equations in field theory (\cite{GIM, Gotay, CCI}).

As a generalization of \cite{tese,FG}, this paper is devoted to the study of
differential forms on a manifold $P$ with local isotropic and decomposable
distributions associated to them, that is, we shall consider a general
$(n+1)$-form $\omega$ on a manifold $P$ that at each point $p$ of $P$ admits
a subspace $L_p$ of $T_pP$ having the following algebraic properties: 
\begin{enumerate}
 \item[(i)] It is \textbf{maximal isotropic} with respect to $\omega_p$: 
       \begin{equation}
        \forall \, v, u \in L_p \qquad \ic_{v \smwedge u}\, \omega_p = 0 \, ,
       \end{equation}
       and $L_p$ is maximal in $T_pP$ with this property.
 \item[(ii)] It is \textbf{decomposable} with respect to $\omega_p$:  $L_p$ has a basis 
       $\{v_1, \ldots, v_m\}$ satisfying 
       \begin{equation}
        \ic_{v_i}\, \omega_p \qquad   
        \text{is a decomposable $n$-form for each $i=1, \ldots, m$.}
       \end{equation}
\end{enumerate}
 The main goal of this article is to show that whenever $L$ is an integrable distribution, 
we can represent $\omega$ locally as a sum
\begin{equation} \label{DarbouxF}
  \omega = \Omega_{\mathcal{F}} + \omega_{\mathcal{F}}
\end{equation}
where $\mathcal{F}$ is a local foliation in $P$ such that $TP = L \oplus T\mathcal{F}$,
$\omega_{\mathcal{F}}$ is the restriction of $\omega$ to $\mathcal{F}$ and 
$\Omega_{\mathcal{F}}$ can be represented by coordinates preserving the decomposition
$TP = L \oplus T\mathcal{F}$ and in which its coefficients are constant, provided it is closed.
As its main application we determine a necessary and sufficient
condition to the existence of 'canonical coordinates' for differential forms. As a special
case, we determine when it is possible to find coordinates like the ones in equation
(\ref{eq:MSPLF1}).

The principal novelties presented here are:
\begin{itemize}
 \item The entire study is made in the general context of a manifold, different from
       (\cite{tese,FG}), where fiber bundles are used from the beginning.
 \item The approach given here to characterize a differential form $\omega$ which has
       canonical coordinates $(x^\mu, q_i, p, p_i^\mu)$ just like (\ref{eq:MSPLF1})
       is extended naturally, and with no extra effort, to include the degenerate case given
       by the restriction of $\omega$ to any submanifold characterized by the relations       
       \begin{equation} \label{eq:constraint}
        dq_{}^i = 0\, , \quad dp\>\!_k^\mu = 0 \, , \qquad \text{for some indexes} 
        \quad i, k, \mu \, .
       \end{equation}
       Under the existence of these constraints $\omega$ is decribed by
       \begin{equation} \label{eq:MSPLF1constraint}
        \omega~=~\sum_{i,\mu \in I} \, dq_{}^i \>\smwedge\, 
          dp\>\!_i^\mu \>\smwedge\, d^{\,n} x_\mu^{} \, - \,
          dp \;\smwedge\, d^{\,n} x~,
       \end{equation}
       where $I$ is the subset of indexes that exclude those appearing in the constraints
       (\ref{eq:constraint}).
 \item The term '$\omega_{\mathcal{F}}$' in the r.h.s. of equation (\ref{DarbouxF}) is new
       and appears as a "horizontal" obstruction to describe $\omega$ with coordinates
       in which it has constant coefficients. It might be possible that, under certain
       conditions, it describes the analog of a background electromagnetic field
       coupled to the classical fields, as in (\cite{Baez}), where we have to add
       a 'horizontal pertubation' $\pi^*\,\omega_F$ to the form $\Omega_0$, 
       the last one being the canonical form (\ref{eq:MSPLF1}), 
       \begin{equation} \label{eq:eletromagnetic}
         \omega =  \Omega_0 + \pi^*\,\omega_F \, ,
       \end{equation}
       to get the right description of a bosonic string in a background electromagnetic 
       field. In this example $\pi$ is a surjective map and $F$ is the electromagnetic 
       field strength.   
\end{itemize}  

We proceed as follows: in section \ref{sec:DEGA} we review some basic algebraic concepts. 
In section \ref{subsec:MIDS} we go further in multilinear algebra and derive the key
results that will be of fundamental importance later. In section \ref{sec:MM} we
start the differential geometric description and prove our main theorem before we show
how to translate our results to fibered manifolds. In the last section we work with 
examples, showing that differential forms with maximal isotropic and decomposable
distributions are quite general, although they do not seem to be generic.


\section{Basic definitions in the Grassmann Algebra} 
  \label{sec:DEGA}
Let $F$ be a finite dimensional vector space and $F^*$ its dual.   
Fixing the notation, if $S \subset F^*$ is any subset, then the subspace defined by
\begin{equation}
 S^\bot  := \{ v \in F \, | \quad 
        \ic_v \, \alpha = 0 \quad \forall \, \alpha \in S \}
\end{equation}
is the annihilator of the set $S$.

 For a $n$-form $\beta$ on $F$,i.e., $\beta \in \bwedge^n \, F^*$, 
we define the contraction map by
\begin{equation}
  \beta^\flat : F  \to \bwedge^{n-1}\,F^*  \qquad \beta^\flat (v) := \ic_v \, \beta \, ,
\end{equation}    
its \textbf{kernel} to be the subspace of $F$ given by 
\begin{equation}
 \ker \beta := \{ v \in F \, | \quad \ic_v \, \beta  = 0 \, \} 
\end{equation}
and its \textbf{support} to be the subspace of $F^*$ given by    
\begin{equation}
  S_\beta  = \bigcap \, \{ \, S \quad | \quad 
         \text{$S$ is a subspace of $F^*$ and}\, \beta \in \bwedge^n \, S \} \, ,
\end{equation}
which is well defined, since the wedge product has the property 
\begin{equation}
  \bwedge^n \, S_1 \cap \bwedge^n \, S_2 = \bwedge^n \, (S_1 \cap S_2)
\end{equation}
for any subspaces $S_1, S_2$ of $F^*$. Hence, $S_\beta$ is the 
smallest subspace of $F^*$ such that $\beta \in \bwedge^n \, S_\beta$.
One can verify that (\cite{tese,FG})
\begin{equation} \label{eq:kerSup}
 S_\beta^\bot = \ker \beta  \quad \text{and} \quad S_\beta = (\ker \beta)^\bot \, .
\end{equation}
Therefore
\begin{equation}
 \dim F = \dim S_\beta + \dim \ker \beta \, .
\end{equation}
When $\beta$ has a trivial kernel we say that it is \textbf{non-degenerate}, 
which is equivalent to say that $S_{\beta}=F^*$.  


Pick a basis $\mathfrak{B}=\{e^1, \ldots , e^{n+N}\}$ of $F^*$, put 
\begin{equation} \label{eq:decomposable2}
 \beta = a_{i_1 \dots i_n} \, e^{i_1} \smwedge \ldots \smwedge e^{i_n} \, , 
\end{equation}
and define $\ell_{\mathfrak{B}}(\beta)$ to be the length of $\beta$ with respect
to $\mathfrak{B}$, that is, 
\begin{equation} \label{eq:lengthbasis}
  \ell_{\mathfrak{B}}(\beta) =  \# \{ \vec{i} \in  \Im^{n+N}_n \, | \,
    \quad a_{i_1 \dots i_n} \neq 0 \} \, ,
\end{equation}
where we use the notation
\begin{equation} \label{indice1}
  \Im^{n+N}_n = \{ \vec{i} := \, 
           (i_1, \ldots, i_n) \, | \quad 1\le i_1 \ldots < i_n \le n+N \}
  \qquad n+N=\dim F \, . 
\end{equation}
The \textbf{length} of $\beta$ is the minimum of the relative
lengths among all possible basis, i.e., 
\begin{equation} \label{eq:length}
 \ell (\beta) = \min_{\mathfrak{B}} \ell_{\mathfrak{B}}(\beta) \, .
\end{equation}
We say that $\beta$ is \textbf{decomposable} if %
\footnote{Trivially, $\ell (\beta)=0$ implies $\beta = 0$.}
$\ell (\beta) \le 1$. If $\ell (\beta)=1$, there is a L.I. set
$\{\alpha^1, \ldots , \alpha^n\}$ such that
\begin{equation} \label{eq:decomposable}
 \beta = \alpha^1 \smwedge \ldots \smwedge \alpha^n
\end{equation}
This set forms a basis for $S_\beta$. Therefore, for any 
$\beta \ne 0$, $\dim S_\beta = n$ if, and only if, $\beta$ is decomposable
(For additional information on decomposable elements see \cite{Marcus}).



\section{Maximal Isotropic Decomposable Subspaces}
\label{subsec:MIDS}

 Let $W$ be a finite dimensional vector space, $L$ a subspace of $W$, 
$\omega$ a $(n+1)$-form on $W$ and $k$ an integer satisfying 
$\, 0 \le k \le n$. The \textbf{$k$-orthogonal complement} of~$L$ in~$W$ 
with respect to $\omega$ is the subspace of $W$ given by
\begin{equation}
 L^{\omega,k}~=~\{ \, v \smin W~|~\mathrm{i}_v^{} \mathrm{i}_{v_1}^{} \ldots
                                 \mathrm{i}_{v_{k}}^{} \omega = 0~~
                     \mbox{for all $\, v_1,\ldots,v_{k} \smin L$} \, \}~.
\end{equation}
In the case $k = 0$ we have $L^{\omega,0}=\ker \omega$. 
The subspace $L$ is said to be, with respect to $\omega$, 
\begin{enumerate}
  \item[(i)] \textbf{$k$-isotropic} when $L \subset L^{\omega,k}$. For $k=1$ we say just
             \textbf{isotropic};
  \item[(ii)] \textbf{strict $k$-isotropic} if it is $k$-isotropic but 
              not ($k - 1$)-isotropic ;  
  \item[(iii)] \textbf{maximal $k$-isotropic} if it is $k$-isotropic and 
               not a proper subspace of another $k$-isotropic subspace; 
\end{enumerate}   
Note that 
\begin{equation} 
 \text{$L$ is maximal $k$-isotropic} 
 \quad \Rightarrow \quad \ker \omega \subset L \, ,
\end{equation}
for, in this case, we have that $\ker \omega + L$ is $k$-isotropic, and so, 
$\ker \omega + L \subset L$ by the maximality condition. 

A vector $v \in W$ is \textbf{decomposable} with respect to $\omega$ if
\begin{equation}
  \ic_v \omega \in \bwedge^n \, W^* \quad 
  \quad \text{is decomposable}. 
\end{equation}
We say that a subspace $L$ is \textbf{decomposable} with respect to $\omega$
if there is a basis for $L$ made of decomposable vectors, with respect to $\omega$.

Any subspace spanned by $\{e_1, \ldots, e_k\}$, with $k \le n$, 
is $k$-isotropic (and strict $k$-isotropic if 
$\ic_{e_1} \ldots \ic_{e_k}\omega \neq 0$). Therefore, maximal $k$-isotropic
subspaces always exist, but are not necessarily decomposable if $n>1$. 
When $k=n=1$ and $\omega$ is non-degenerate,
maximal isotropic subspaces are called Lagrangian and they always admit a complementary 
maximal isotropic subspace, which is lagrangian as well. Extending this result for forms
of any degree:


\begin{thm}
 \label{Thm:DME} \mbox{} \\[1mm]
  Let $L$ and $V$ be subspaces of $W$ such that $L \subset V$, $L$ is maximal isotropic decomposable and 
  $V$ is $r$-isotropic with respect to a ($n+1$)-form $\omega$ on $W$.
  There is a $n$-isotropic subspace $F$ such that 
  \begin{equation} \label{Thm:Dec}
    W = L \oplus F \qquad \quad (F \cap V) \, \oplus \, L = V   
  \end{equation}
  and $F \cap V$ is ($r-1$)-isotropic.
\end{thm}

\proof

Assume that $\omega$ is non-degenerate, i.e., $\ker \omega=0$.
This proof is by induction on $m+1=\dim L$.
Pick a decomposable basis 
$\{ v_0, \ldots, v_m \}$ for $L$ and a $1$-form
$\alpha_0 \in W^*$ such that $\alpha_0(v_i)=\delta_i^0$.
There are $1$-forms $u^1, \ldots, u^n \in L^\bot$ such that 
\[
  \ic_{v_0}\omega = u^1 \smwedge \ldots \smwedge u^n \, ,
\]
with at least ($n-r+1$) of them in $V^\bot$, 
since $L \subset V$ and $V$ is $r$-isotropic.
Define $L_1 = L \cap \ker \alpha_0$, which is generated by $\{ v_1, \ldots, v_m \}$, 
and the $(n+1)$-form $\omega_1$ 
\begin{equation} \label{eq:dec}
 \omega = \omega_1 + \, \alpha_0 \, \smwedge u^1 \smwedge \ldots \smwedge u^n .
\end{equation}

Take any subspace %
\footnote{Notation: ``$\left\langle u_1, \ldots, u_k \right\rangle$''
is the subspace spanned by the L.I. set $\{ u_1, \ldots, u_k \}$. } 
$F_1= \left\langle u_1, \ldots, u_k \right\rangle \subset \ker \alpha_0$,
such that 
\begin{enumerate}
 \item[(i)] $\ker \omega_1 = \left\langle u_1, \ldots, u_s \right\rangle 
             \oplus \left\langle v_o \right\rangle$ , $s \le k$ ;
 \item[(ii)] $\ic_{u_1 \smwedge \ldots \smwedge u_k}u^1 \smwedge \ldots \smwedge u^n \ne 0$ ;
 \item[(iii)] If $n \ne s$ then $\ic_{u_{s+1} \smwedge \ldots \smwedge u_k}\omega_1 \notin \bwedge \, L^\bot$ ;
\end{enumerate}
Such a nontrivial $(k \ge 1)$ subspace $F_1$ always exists. For $k \le n$ 
maximal satisfying such properties, we shall prove that $k=n$. If $n=s$, it is trivial. 
If $s \le k < n$, there is a vector 
$u_{k +1} \in \ker (\alpha_0 \wedge u^1 \smwedge \ldots \smwedge u^s)$ such that
\[
  \ic_{u_{s+1} \smwedge \ldots \smwedge u_{k+1}}\omega_1 \notin \bwedge \, L^\bot
\]
by the property (iii) and the fact that $L$ is isotropic. 
Applying the same reasoning recursively, we get $k=n$. Therefore, 
\begin{equation} \label{eq:F_0}
 \ic_{u}\omega_1 \in \bwedge^n L_1^\bot \quad \Longrightarrow 
 \quad \ic_{u}\omega_1 = 0 \, .
\end{equation}
for each $u \in F_1$. 

Define 
$W_1 = \ker (\alpha_0 \wedge u^1 \smwedge \ldots \smwedge u^s)$, where we 
are assuming that $F_1= \left\langle u_1, \ldots, u_n \right\rangle$ has
the first $s$ vectors, together with $v_0$, composing the kernel of $\omega_1$
and $u^j{(u_i)=\delta^j_i}$.
Note that $\omega_1$ is non-degenerate in $W_1$ and $L_1 \subset W_1$
is isotropic and decomposable w.r.t. $\omega_1$, since
$L_1=\left\langle v_1, \ldots, v_m \right\rangle$.
To see that $L_1$ is maximal isotropic w.r.t. $\omega_1$ in $W_1$,
let $e = e_1 + u \in W_1$ be written according to the decomposition
\[
  W_1= E\oplus U
       \quad U =\left\langle u_{s+1}, \ldots, u_n \right\rangle 
       \quad E = \ker (\alpha_0 \wedge u^1 \smwedge \ldots \smwedge u^n)
\]
with $e_1 \notin L_1$. By the maximal $1$-isotropy of $L$ w.r.t. 
$\omega$ in $W$, there are vectors $e_2, \ldots e_n \in W_1$ and $v \in L_1$
such that 
\[
  \omega(v,e_1, \ldots, e_n)=\omega_1(v,e_1, \ldots, e_n)=1 
\] 
Furthermore, noticing that the annihilator $L^\bot$ of $L$ is generated
by the $1$-forms $\ic_{v' \smwedge e_2' \smwedge \ldots \smwedge e_n'}\,\omega_1$,
with $e_2', \ldots e_n' \in W/L$ and $v' \in L$, we can assume that
\[
  \omega(v,u,e_2, \ldots, e_n)=\omega_1(v,u,e_2, \ldots, e_n)=0 \, . 
\] 
Therefore,
\begin{equation} \label{eq:FRE}
  \ic_{e_1+u} \omega_1 \notin \bwedge^n L_1^\bot. 
\end{equation}
So, 
if $e=e_1 + u \in W_1$ is such that
\[
 \ic_{e \smwedge v}\omega_1 = 0
\]
for all $v \in L_1$, then $e_1 \in L_1$, and by (\ref{eq:F_0}), $u=0$. 
In other words, $L_1$ is maximal
isotropic w.r.t. $\omega_1$ in $W_1$, implying that it is
maximal isotropic decomposable w.r.t $\omega_1$. By the induction hypothesis, let $F' \subset W_1$
be $n$-isotropic w.r.t. $\omega_1$ and complementary to $L_1$ in
$W_1$. Then $F= \left\langle u_1, \ldots, u_s \right\rangle \oplus F'$
is $n$-isotropic w.r.t. $\omega$ and complementary to $L$ in $W$. Furthermore,
$F \cap V \cap W_1$ is ($r-1$)-isotropic, and so is  $F \cap V$.

\qed

An useful lemma that also helps understanding the content of a maximal isotropic decomposable
subspace is given bellow.


\begin{lem} \label{lem:useful}  
Let $L$ be a maximal isotropic decomposable subspace with respect to a 
($n+1$)-form $\omega$ on $W$, and  $f_1, \ldots, f_{n}$ vectors in 
$W$ such that
\[
  \ic_{f_{1} \smwedge \dots \smwedge f_{n}}\omega \notin L^\bot 
\]
Then there are $1$-forms $f^1, \ldots, f^{n}$ such that $f^j(f_i) = \delta^j_i$ and 
$f^1 \smwedge \ldots \smwedge f^{n} \in \omega^\flat (L)$. \\
\end{lem}
\proof 
Since $\ic_{f_{1} \smwedge \dots \smwedge f_{n}}\omega \notin L^\bot$,
there is a $\omega$-decomposable vector $v \in L$ such that 
$\omega(v, f_{1}, \dots , f_{n}) = 1$. Since $\ic_v \, \omega$ is decomposable,
its kernel has codimension $n$ in $W$, and therefore we can write
\[
  W = \left\langle f_{1}, \dots , f_{n} \right\rangle \oplus 
       \ker \ic_v \, \omega \, . 
\]
This implies that for $f^1, \ldots, f^{n} \in (\ker \ic_v \, \omega)^\bot$ such that 
$f^j(f_i) = \delta^j_i$ we have $\ic_v \, \omega = f^{1} \smwedge \dots \smwedge f^{n}$.
\qed

Let  $L$  be a maximal isotropic decomposable subspace and $F$ a complementary $n$-isotropic 
subspace in $W$, both with respect to $\omega$. Pick any basis 
$\mathfrak{B}_{F}  = \left\{ f^{}_1, \ldots, f^{}_{N+n} \right\}$ of $F$
and define the number of non-vanishing indexes of 
$\mathfrak{B}_{F}$ with  respect to $\omega$ :
\begin{equation} \label{indiceF2}
  \mathfrak{N}(\mathfrak{B}_F) = \# \, \Im^{\mathfrak{B}_F}
  \qquad \text{where} \quad 
  \Im^{\mathfrak{B}_F} =  \{ \, \vec{i} \in  \Im^{n+N}_n \, | \,
     \ic_{f_{i_1} \smwedge \dots \smwedge f_{i_n}}\omega \neq 0 \}
\end{equation}
and $\Im^{n+N}_n$ is given by (\ref{indice1}). 
Let $\mathfrak{B} = \left\{ e_1, \ldots , e_m, f_1, \ldots, f_{N+n} \right\}$ 
be a basis of $W$ adapted to the decomposition (\ref{Thm:Dec}), that is,
$\left\{ e_1, \ldots , e_m \right\}$ is a basis of $L$ and  
$\mathfrak{B}_F \subset \mathfrak{B}$, and let
$\mathfrak{B}^* = \left\{ e^1, \ldots , e^m, f^1, \ldots, f^{N+n} \right\}$ 
be its dual. Hence
\begin{equation} \label{eq:RepresGeralCoord}
  \omega = \sum_{\vec{i} \in \Im^{\mathfrak{B}_F}} 
     \alpha_{\vec{i}} \, \smwedge f^{i_1} \smwedge \ldots \smwedge f^{i_k} \, ,
     \qquad \quad  \alpha_{\vec{i}} \in F^\bot \cong L^*\, ,
\end{equation} 
where for each $\vec{i} \in \Im^{\mathfrak{B}_F}$, $\alpha_{\vec{i}}$ is a non-trivial
linear combination of the $1$-forms $e^1, \ldots , e^m$. %
Therefore,  
\begin{equation} \label{eq:RepresGeralCoord2}
 \mathfrak{N}(\mathfrak{B}_F) = \ell_{\mathfrak{B}} (\omega)
 \qquad \text{, which implies} \qquad
 \mathfrak{N}(\mathfrak{B}_F) \ge \ell (\omega) \, .
\end{equation} 
where $\ell_{\mathfrak{B}} (\omega)$ is the length of
$\omega$ with respect to $\mathfrak{B}$, defined by formula (\ref{eq:lengthbasis}).
Taking the minimum among all possible basis of $F$,
we can define the number %
\footnote{The conclusion ``$\mathfrak{N}_L \ge 0$'' follows from formula 
(\ref{eq:RepresGeralCoord}).}  
\begin{equation} \label{indiceF4}
  \mathfrak{N}_L = \min_{\mathfrak{B}_F} \,
  \, \mathfrak{N}(\mathfrak{B}_F) - \dim (L / \ker \omega ) 
  \ge 0\, . 
\end{equation}
It is not difficult to check that $\mathfrak{N}_L$ does not depend 
on the choice of the complementary $n$-isotropic subspace.


\begin{thm}
\label{thm:TLCC}  \mbox{} \\[1mm]
Let $L$, $F$ and $V$ be subspaces of $W$ satisfying the relations in 
theorem \ref{Thm:DME} with respect to $\omega \in \bwedge^n \, W^*$, 
and suppose that $\mathfrak{N}_L = 0 $, which is equivalent to
\begin{equation} \label{eq:thm:TLCC1}
 \dim (L/ \ker \omega) = \ell (\omega) \, .
\end{equation}
There is  a basis $\mathfrak{B}_F$ for $F$ adapted  to $V$, with 
$\mathfrak{B}_F^*=\left\{ e^1, \ldots, e^{N+n} \right\}$ its unique dual 
in $L^{\bot}$, and a L.I. set $\left\{ \hat{e}_{\vec{i}_1}, \ldots, \hat{e}_{\vec{i}_m} \right\}$
in $F^{\bot}$ such that
\begin{equation} \label{eq:FML}
  \omega~= \sum_{\vec{i} \in \Im^{\mathfrak{B}_F}} \, 
           \hat{e}\>\!_{\vec{i}}^{} \,
           \smwedge\, e^{i_1} \smwedge \ldots \smwedge\, e^{i_n}
\end{equation} 
where $\Im^{\mathfrak{B}_F}$ is given by equation (\ref{indiceF2}). Each term 
$e^{i_1} \smwedge \ldots \smwedge\, e^{i_n}$ in equation (\ref{eq:FML}) vanishes
when contracted with $r$ vectors in $V$. 
\end{thm}

\proof
Take a basis $\mathfrak{B}_F$ for $F$ with $\#\Im^{\mathfrak{B}_F} = \dim L - \dim \ker \omega$,
which exists since $\mathfrak{N}_L=0$. Define  
\[
  \hat{e}\>\!_{\vec{i}} = \ic_{e_{i_1} \smwedge \ldots \smwedge\, e_{i_n}}\omega
\]
Using the hypothesis on the dimension of $L/ \ker \omega$, the $n$-isotropy of $F$
and the the $1$-isotropy of $L$, it is easy to check that they form a basis for
$(L/ \ker \omega )^*$ and that equation (\ref{eq:FML}) holds.
Relation (\ref{eq:thm:TLCC1}) follows 
from the formula (\ref{eq:FML}).
\qed



\section{Maximal Isotropic Decomposable Distributions}
\label{sec:MM}



\subsection{Maximal Isotropic Decomposable Distributions and Flatness}
\label{subsec:LMDFDarboux}

Let $\omega \in \Omega^{n+1}(P)$ be a ($n+1$)-differential form on a manifold $P$. 
We say that it is \textbf{flat} if, around
each point of $P$, there is a coordinate representation in which
it has constant coefficients, that is, 
\begin{equation}
 \omega =  \omega_{i_1 \ldots i_{n+1}} 
 \, dx^{i_1} \smwedge \ldots \smwedge \, dx^{i_{n+1}} 
 \qquad \text{and}  \qquad
 {\partial \over \partial x^i} \, \omega_{i_1 \ldots i_{n+1}} \equiv 0 
 \qquad \text{, for each $i, i_1, \ldots, i_{n+1}$.} 
\end{equation} 

A distribution $L$ on $P$ is \textbf{maximal isotropic decomposable} with respect to $\omega$ 
if it is pointwise maximal isotropic decomposable with respect to $\omega$.

\begin{thm}[Principal Part of $\omega$]\label{ex:ASF}
 \mbox{} \\[1mm] 
If  $L$ is a maximal isotropic decomposable distribution with respect to $\omega$, then 
for any differential form $\delta\omega \in \Omega^{n+1}(P)$, we have
\begin{equation} \label{eq:ExempAddingForm}
 \delta\omega^\flat(L) = 0 \qquad \Rightarrow \qquad
 \text{$L$ is maximal isotropic decomposable with respect to $\omega + \delta\omega$} \, .
\end{equation}   
Therefore, $L$ defines a class of forms admitting it as a 
maximal isotropic decomposable distribution,
\begin{equation}
 [\omega]_L = \{\, \omega + \delta \omega \quad | \quad 
               \delta\omega^\flat(L) = 0  \, \} \, 
\end{equation} 
which we will call the \textbf{principal part} of $\omega$ with respect to $L$.
\end{thm}
\proof
To check this, note that $L$ is isotropic and decomposable with respect to
$\omega + \delta\omega$, and if $u$ is a vector fields on $P$ such that
$\, \ic_{u \smwedge v} (\omega + \delta\omega) = 0 \,$ 
for every $v \in L$,we have
\begin{equation} 
 \ic_{u \smwedge v} \omega = \ic_{u \smwedge v} (\omega + \delta\omega) =  0\, ,
\end{equation}   
for every $v \in L$, implying that $u \in L$, since $L$ is maximal isotropic 
with respect to $\omega$. Then, $L$ is also maximal isotropic 
with respect to $\omega + \delta\omega$. 
\qed


\begin{thm}[Flatness] \label{Thm:MFlat}
  \mbox{} \\[1mm] 
 Let $\omega \in \Omega^{n+1}(P)$ 
 satisfies the regular condition of constant dimension of
 the kernel distribution. 
 If $L$ is a maximal isotropic decomposable distribution, then for
 each foliation $\mathcal{F}$ such that $TP=T\mathcal{F} \oplus L$, 
 we have that
 \begin{equation}
  \Omega_{\mathcal{F}}=\omega - \omega_{\mathcal{F}} \in [\omega]_L \, ,
 \end{equation}
 where $\omega_{\mathcal{F}}$ is the restriction of $\omega$ to $\mathcal{F}$,
 defines a representation of the principal part of $\omega$ and admits $T\mathcal{F}$
 as a $n$-isotropic distribution. If $L$ is integrable, then 
 \begin{center}
  $\Omega_{\mathcal{F}}$ is closed if, and only if, it is flat.
 \end{center}
 Moreover, a coordinate system for $P$ on which  
 $\Omega_{\mathcal{F}}$ has constant coefficients can be chosen to be adapted %
 \footnote{That is, cordinates $(x,y)$ such that $T\mathcal{F}$ 
 is given by $dy=0$ and $L$ by $dx=0$.}
 to the decomposition $TP=T{\mathcal{F}} \oplus L $. 
\end{thm}

\proof

Since the character of this theorem is local, we will avoid 
to use the label ``local'', keeping in mind that there is no need 
for global constructions here. 
Hence, we will assume that $P=X \oplus L$,
is a vector space where $L$ is a subspace identified with the 
maximal isotropic decomposable distribution and $X$ a subspace such that
$X \times \{y \}$ are the leaves 
of the foliation $\mathcal{F}$, which is always possible since $L$ and
$T\mathcal{F}$ are simultaneously integrable 
(see appendix \ref{appendix:Integrability}). 
Moreover, since the distribution $L$ is maximal isotropic decomposable 
with respect to $\Omega_{\mathcal{F}}$ (see lemma \ref{ex:ASF}), 
we can prove the theorem for the case $\omega=\Omega_{\mathcal{F}}$,\
that is, $\omega_{\mathcal{F}}=0$.
Therefore the foliation $\mathcal{F}$ will be, just in this case,
$n$-isotropic with respect to $\omega$.

Let $\omega_0^{}$ be the constant form obtained by spreading
$\omega(p_0)$, the value of $\omega$ at the origin $p_0=(0,0)$, all over the 
vector space $X \oplus L$. Clearly, the distribution $L$ is maximal isotropic decomposable 
with respect to $\omega_0$.


\begin{lem} \label{lem3:ProofDarboux}
 In a small neighborhood of the origin, 
 $\qquad \omega_0^\flat (L) = \omega^\flat (L)$. 
\end{lem}
\proof
The vector fields $\{ f_\mu:= {\partial \over \partial x^\mu} \}$ 
form a basis for the $n$-isotropic distribution $T\mathcal{F}$. Define the $1$-forms
\begin{equation}
 \alpha_{\mu_1 \ldots \mu_n}:= 
 \ic_{f_{\mu_1} \smwedge \ldots  \smwedge f_{\mu_n}} \omega 
\qquad \quad
 \alpha^0_{\mu_1 \ldots \mu_n}:=  \alpha_{\mu_1 \ldots \mu_n} (p_0) =
 \ic_{f_{\mu_1} \smwedge \ldots  \smwedge f_{\mu_n}} \omega_0 
\end{equation}
By continuity, there is a small neighborhood of the origin 
where we have that $\alpha^0_{\mu_1 \ldots \mu_n} \ne 0$ implies 
$\alpha_{\mu_1 \ldots \mu_n} \ne 0$, that is, for each point $p$
in this neighborhood
\begin{equation}
 (\ic_v \, \omega_0)(v_1 \smwedge \ldots  \smwedge v_n)  \ne 0
 \quad \Rightarrow \quad
 (\ic_v \, \omega)(v_1 \smwedge \ldots  \smwedge v_n) \ne 0
\end{equation}
for any vectors $v_1, \ldots , v_n \in T_p\mathcal{F}$ and
$v \in L_p$. This implies that the annihilators satisfy
$(\omega^\flat (L_p))^\bot \subset (\omega_0^\flat (L_p))^\bot$,
therefore
$\omega_0^\flat (L_p) \subset \omega^\flat (L_p)$.
Since $\dim \ker \omega$ is constant by 
hypothesis and $L$ is integrable, 
$\dim \omega^\flat (L) = \dim L/ \ker \omega$ must be constant,
and so, there must exist a small neighborhood where
$\omega_0^\flat (L) = \omega^\flat (L)$. 

\qed


\begin{lem} \label{lem2:ProofDarboux}
 $\qquad \exists \quad \theta \, , \theta_0 \in \omega^\flat (L) \quad : 
 \quad \omega = d\theta \quad \omega_0 = d\theta_0$
\end{lem}
\proof
Analogously to the ``canonical'' proof of Poincare lemma, define the ``$L$-contraction'' 
$\Phi_t(x,y)=(x,ty)$, for each $\, t \smin \mathbb{R}$. Denote its time dependent 
vector field by $\xi_t$, that is, 
\begin{equation}
 \xi_t \bigl( \Phi_t(x,y) \bigr)~=~\frac{d}{ds} \, \Phi_s(x,y) \, \bigg|_{s=t}
   \qquad  \xi_t(x,y)~=~t^{-1}(0,y)  \, . 
\end{equation}
Although the vector field $\xi_t$ have a singularity at $t=0$,
the forms 
\begin{equation}
 \theta_t=(\Phi_t)^*(\mathrm{i}_{\xi_t}^{} \omega) 
\end{equation}
are well defined for every $t \in \R$. Since 
$\xi_t \in L$, we have $\mathrm{i}_{\xi^L_t}\omega \in \omega^\flat (L)$,
that is, $\theta_t \in \omega^\flat (L)$. Defining the $n$-form
\begin{equation}
 \theta~  =~\int_0^1 dt~\theta_t \, , 
\end{equation}
it is clear that $\theta \in \omega^\flat (L)$.
Moreover, since $\, d\omega = 0 \,$, $\Phi_1=Id_P$, and $(\Phi_0)^*(\omega)=0$,
the last relation following by the $n$-isotropy of $T\mathcal{F}$, we have 
 \begin{eqnarray*}
  d\theta \!\!
   =\!\! \int_0^1 dt~
          (\Phi_t)^*d(\mathrm{i}_{\xi_t}^{} \omega)  
   =~     \int_0^1 dt~ (\Phi_t)^*(L_{\xi_t}^{} \omega) 
   = \!\! \int_0^1 dt~
          \frac{d}{dt} \bigl( (\Phi_t)^* \omega \bigr) \, 
   =~ \omega \, .    
 \end{eqnarray*}
We can apply the theorem for $\omega_0$, and since 
$\omega_0^\flat (L)=\omega^\flat (L)$, prove the assertion of the lemma.
\edem

Consider the family of $(n+1)$-forms given by 
$\, \omega_t^{} =  \omega_0^{} + t(\omega-\omega_0^{})$, for every 
$\, t \smin \mathbb{R}$. 
By lemma \ref{lem3:ProofDarboux} we have
$\omega_t^\flat (L) \subset \omega^\flat (L)$, and the continuity 
of the argument implies that 
there is an open neighborhood of the origin $p_0=(0,0)$ where, 
for all $t$ satisfying $\, 0 \le t \le 1 \,$ and all points $p$ in it,
$\omega_t^{}(p)$ satisfies 
\begin{equation} \label{eq:EqualContraction}  
 \omega_t^\flat (L_p) = \omega^\flat (L_p) \, . 
\end{equation}
By lemma \ref{lem2:ProofDarboux}, we have $n$-forms 
$\theta, \theta_0 \in \omega^\flat (L)$ satisfying $\omega = d\theta$ and 
$\omega_0 = d\theta_0$. Therefore, 
\begin{equation}
 \alpha := \theta_0 - \theta  \quad
 \text{satisfies} \quad 
 \alpha \in \omega_t^\flat (L) \quad \text{and} \quad
 \quad d\alpha = \omega_0^{} - \omega \, ,
\end{equation}
for every $\, 0 \le t \le 1 \,$ and every point in this ``small'' neighborhood 
of the origin where the relation (\ref{eq:EqualContraction}) holds.
This implies that we can pick a time dependent vector field $X_t^{} \in L$ 
defined on this neighborhood satisfying
\[
  \mathrm{i}_{X_t}^{} \omega_t^{}~=~\alpha~.
\]
Let $\, \Phi^X_t \equiv \Phi^X_t((0,0))$ be its flux beginning at the origin,
which is well defined for $\, 0 \le t \le 1$, in some open neighborhood of the 
origin. Then it follows that
 \[
 \begin{array}{rcl}
  {\displaystyle \frac{d}{ds} \, \bigg|_{s=t} (\Phi^X_s)^* \omega_s^{}}
  &=& {\displaystyle (\Phi^X_t)^* \bp \frac{d}{ds} \bigg|_{s=t} \omega_s^{} \ep +
                     \frac{d}{ds} \bigg|_{s=t} (\Phi^X_s)^* \omega_t^{}} \\[5mm]
  &=& (\Phi^X_t)^* \bigl( \omega - \omega_0^{} + L_{X_t}^{} \omega_t^{} \bigr) \\[3mm]
  &=& (\Phi^X_t)^* \bigl( \omega - \omega_0^{} +
                   d(\mathrm{i}_{X_t}^{} \omega_t^{}) \bigr) \\[3mm]
  &=& (\Phi^X_t)^* \bigl( \omega - \omega_0^{} + d\alpha \bigr) \\[3mm]
  &=& 0
 \end{array}
 \]
Therefore, $\Phi^X_1$ is the desired coordinate transformation, since
$\, (\Phi^X_1)^* \omega = (\Phi^X_1)^* \omega_1^{} = (\Phi^X_0)^*\omega_0^{} = \omega_0^{}$.
Since $X_t$ is in $L$, its flux acts as the identity in the leaves of $\mathcal{F}$.
\qed



\subsection{Maximal Isotropic Decomposable Distributions on Fibered Manifolds}
\label{subsec:LMDFFibred}

 Let $P$ be a fibered manifold over $M$, that is, a surjective 
submersion $P \stackrel{\pi}{\to} M$. 
Denote its vertical distribution by $V := \ker T\pi$. 

 When $\mathfrak{N}_L \equiv 0$, we can apply theorems  \ref{thm:TLCC} and 
\ref{Thm:MFlat} to find ``canonical coordinates'' for the principal
part of $\omega$. Before stating this theorem, let's first fix the 
notation: for each $0 \le s \le n$ 
\begin{equation} \label{IndAdap}
  \Im^s_n(\pi,L) = \{ \vec{i}_s \, \times \vec{\mu}_s := \, 
           (i_1, \ldots, i_s, \mu_1, \ldots, \mu_{n-s})\}_{\quad
           1\le i_1 \ldots < i_s \le N'}^{\quad 1 \le \mu_1 \ldots < \mu_{n-s} \le n'}
\end{equation}
where $N'= \dim (V/L)$ and $n'=\dim M$. In important applications we have
$n=n'$, but in general it does not obey this equality.


\begin{cor}[Canonical Coordinates] \label{Thm:MCC}
  \mbox{} \\[1mm] 
 Let $\omega \in \Omega^{n+1}(P)$ be a nondegenerate %
 \footnote{This could be changed by the regular condition of constant 
 dimension of the kernel distribution. In this case, we should have 
 $\ell (\omega) + \dim \ker \omega = \dim L \,$ .} 
 form on $P$ and $L$ an integrable maximal isotropic decomposable distribution 
 such that $\mathfrak{N}_L\equiv 0$, that is, the length of $\omega$ satisfies 
 the  relation
 \begin{equation}
  \ell (\omega) = \dim L \, .
 \end{equation} 
 Let $P \stackrel{\pi}{\to} M$ be any fibered manifold, with vertical 
 distribution $V$, and $\mathcal{F}$ any foliation such that 
 \begin{equation} \label{Thm:Dec}
    TP = L \oplus T\mathcal{F} \qquad \quad V=L\oplus (T\mathcal{F}\cap V) \, .   
 \end{equation}
 If $\Omega_{\mathcal{F}}$, 
 the representation of the principal part of $\omega$ given by $\mathcal{F}$, 
 is closed and admits the vertical bundle $V$ as an $r$-isotropic distribution, 
 then, in a small neighborhood of each point of $P$, there are coordinates 
 $(p_j, q^i ,\, x^\mu)\,$ such that
 \begin{equation} \label{eq:FMLBDG2}
  \Omega_{\mathcal{F}}=\omega - \omega_{\mathcal{F}}~= 
  \sum_{s=0}^{r-1} \, \sum_{\vec{i}_s\times\vec{\mu}_s \in \Im^s_\omega} \, 
           dp\>\!_{\vec{i}_s;\vec{\mu}_s}^{} \,\smwedge\,
           dq^{\vec{i}_s} \,\smwedge\, dx^{\vec{\mu}_s} \, ,
 \end{equation}
 where the distribution $T\mathcal{F}$ is given by $dp = 0$, $V$ by $dx=0$ and
 $L$ by $dq=dx=0$. Here we are using the notations: 
 $\omega_{\mathcal{F}}$ is the restriction of $\omega$ to $\mathcal{F}$,
 for each $0 \le s \le r-1$, $\Im^s_\omega$ is a fixed subset 
 of $\Im^s_n(\pi,L)$,
 \begin{equation} 
   dq^{\vec{i}_s} := dq^{i_1} \,\smwedge \ldots \smwedge\, dq^{i_s} \qquad
   dx^{\vec{\mu}_s} := dx^{\mu_1} \,\smwedge \ldots \smwedge\, dx^{\mu_{n-s}} 
   \qquad \text{and} \qquad
   dp\>\!_{\vec{i}_s;\vec{\mu}_s} = dp_{I\,(\vec{i}_s;\vec{\mu}_s)} \, ,
 \end{equation}
 with $I: \Im^0_\omega \cup \ldots \cup \Im^{r-1}_\omega \to \mathbb{N}$ an injection.
\end{cor}

 Now we shall single out the right parameters to make formula (\ref{eq:FMLBDG2}) 
turn into formula (\ref{eq:MSPLF1}). First of all, we must have
\begin{equation} \label{eq:ConditionsCanonico}
  \omega \in \Omega^{n+1}(P) \qquad n= \dim M 
  \qquad N= \dim (V/L) \qquad r=2 \qquad \dim L = Nn + 1.
\end{equation} 
Under these conditions, if we just assume that $L$ is isotropic, we can show that
it is decomposable and maximal isotropic, for the map $\Omega_{\mathcal{F}}^\flat$ takes
$L$ onto the decomposable subspace $\bwedge^n_1 L^\bot$ of $n$-forms  
which vanishes whenever contracted with two vectors in $V$ or any vector in $L$, 
arriving to the conditions
presented in \cite{tese,FG}. According to them, in this case we have: 
\begin{center}
  for $n \ge 1$, $\,L$ is the unique maximal isotropic decomposable distribution in $V$
\end{center} 
and,
\begin{center}
  if $n \ge 2$ and $\Omega_{\mathcal{F}}$ is closed, then $L$ is integrable.
\end{center} 
Moreover, applying corolary \ref{Thm:MCC} with the conditions 
(\ref{eq:ConditionsCanonico}), we can find coordinates 
$(p, p^\nu_j, q^i ,\, x^\mu)\,$ such that $\, i,j = 1, \ldots, N$, 
$\,\mu, \nu = 1, \ldots, n$ and formula (\ref{eq:FMLBDG2}) reads  
\begin{equation} \label{eq:MSPLF1-2}
 \Omega_{\mathcal{F}}~=~ \, dp \;\smwedge\, d^{\,n} x~
        + dp\>\!_i^\mu \>\smwedge\, dq_{}^i \>\smwedge\, d^{\,n} x_\mu^{} \,
\end{equation}
just like the equation (\ref{eq:MSPLF1}). If instead of 
$\Omega_{\mathcal{F}}$, $\omega$ itself
satisfies the hypothesis above, we can choose $\mathcal{F}$, at least locally,
such that $\omega_{\mathcal{F}}=0$, that is, $\Omega_{\mathcal{F}}=\omega$.


\section{Examples}
\label{sec:E&A}


\begin{exem}[Decomposable Forms] \label{Exem:DecomposableForms}
 \mbox{} \\[1mm] 
 If  $\omega \in \Omega^{n+1}(P)$ is a decomposable form on a manifold $P$, a 
local distribution $L$ is maximal isotropic decomposable if, and only if, at each point where
$\omega \ne 0$ 
\begin{equation} \label{eq:ExempDecompMult}
 L = \ker \, \omega \oplus L_0 \qquad  \dim L_0 =1 \, .
\end{equation}   
Therefore, every decomposable form admits a local maximal isotropic decomposable distribution. 
This includes all ($n+1$)-forms if $\dim P \le n+2$. For instance, densities
and volume forms on $P$. 
\end{exem}

\begin{exem}[Product of Manifolds] \label{ex:ProductManifolds}
 \mbox{} \\[1mm] 
Given two manifolds $P_1$ and $P_2$ together with the ($n+1$)-forms 
$\omega_1$ and $\omega_2$ such that they admit (local) 
maximal isotropic decomposable distributions $L_1$ and $L_2$, respectively. Then 
the manifold $P_1 \times P_2$ with the ($n+1$)-form 
$\omega_1 \oplus \omega_2$ admits the (local) maximal isotropic decomposable distribution 
$L_1 \oplus L_2$. 

\end{exem}

\begin{exem}[Isotropic Distribution of Maximal Dimension] \label{Exem:IDMDimension}
 \mbox{} \\[1mm] 
 Every distribution $L$ on a manifold $P$ which is isotropic
with respect to $\omega \in \Omega^{n+1}(P)$, that is, 
$\omega^\flat(L) \subset \bwedge^n \,L^\bot$, clearly satisfies the dimension constraint
\begin{equation} \label{eq:MULTD1} 
   \dim \, L - \dim \, \ker \, \omega = \dim \omega^\flat(L)
   ~\le~ \dim \bwedge^n \,L^\bot \,
          =\, {N+n  \choose n} ~.
\end{equation}
where $N +n = \dim P - \dim L$.  Furthermore, there is the equivalence
\begin{equation} \label{eq:canFormDefinMaxDimension}
 \quad \dim L = \dim \, \ker \, \omega + {N+n  \choose n} 
 \quad \iff \quad
 \omega^\flat (L) = \bwedge^n \, L^{\bot} \, .
\end{equation}   
Therefore, if $L$ has the maximal dimension allowed to a isotropic distribution,
it is maximal isotropic decomposable and $\mathfrak{N}_L=0$. In \cite{Ma2, FG} it is 
proved that such a distribution is unique, and if it is integrable 
and $\omega$ is closed, then there are local
coordinates such that formula (\ref{eq:FMLBDG2}) becomes
\begin{equation} \label{eq:CanCoordTeseMaxDim}
  \omega ~= \sum_{\vec{i} \in \Im^{n+N}_n} \, 
           dp\>\!_{\vec{i}}^{} \,\smwedge\,
           dq^{\vec{i}} \, 
         ~= \sum_{\vec{i} \in \Im^{n+N}_n} \, 
           dp\>\!_{\vec{i}}^{} \,\smwedge\,
           dq^{i_1} \,\smwedge \ldots \smwedge\, dq^{i_n}
\end{equation}
\end{exem}

\begin{exem}[Canonical Examples] \label{ExemCan}
 \mbox{} \\[1mm]  

Here we summarize part of the work done in \cite{tese,FG}.
Let $P \stackrel{\pi}{\longrightarrow} M$ be a fibred manifold
such that its vertical distribution $V$ is $r$-isotropic with respect to 
$\omega \in \Omega^{n+1}(P)$. 
Suppose we have a distribution $L \subset V$ satisfying the relation
\begin{equation} \label{eq:canFormDefinExamp}
 \omega^\flat (L) = \bwedge^n_r \, L^{\bot} \, .
\end{equation}   
This is equivalent to say that $L$ is isotropic and
\begin{equation*} \label{eq:MULTD1EXAMP} 
   \dim \, L~=~\dim \, \ker \, \omega \, + \,
         \sum_{s=0}^{r-1} {N'  \choose s} {n' \choose n-s}~.
\end{equation*}
where ${p \choose q} = 0 \;$ if $\, q > p \,$, $N' = \dim V/L$ and $n' = \dim M$.
Such a distribution is maximal isotropic decomposable with $\mathfrak{N}_L \equiv 0$. 
Furthermore, 
\begin{equation*}
 \text{$L$ is unique in $V$ if} \qquad {n' \choose n+1-r} \ge 2 \, 
\end{equation*}
and for $d\omega =0$, 
\begin{equation*}
 \text{$L$ is involutive if} \qquad {n' \choose n+1-r} \ge 3 \, .
\end{equation*}
In this case it is possible to find a foliation $\mathcal{F}$ such
that $\omega_{\mathcal{F}}=0$ and coordinates such that 
formula (\ref{eq:FMLBDG2}) becomes
\begin{equation*} 
  \omega~= \sum_{s=0}^{r-1} \, 
        \sum_{\vec{i}_s\times\vec{\mu}_s \in \Im^s_n(\pi,L)} \, 
           dp\>\!_{\vec{i}_s;\vec{\mu}_s}^{} \,\smwedge\,
           dq^{\vec{i}_s} \,\smwedge\, dx^{\vec{\mu}_s} \, .
\end{equation*}
Note that the difference between this and formula 
(\ref{eq:FMLBDG2}) is the set of indexes in which they are summed.

\end{exem}

\begin{exem}[A Maximal Isotropic Decomposable Distribution with $\mathfrak{N}_L=1$] 
\label{Exem:N=1} \mbox{} \\[1mm]  
Let $\omega$ be the $4$-form in $\R^{11}$ given by
\begin{equation*} 
  \omega~= dp_1 \,\smwedge \, dq^1 \,\smwedge\, dx_1^{1} \,\smwedge\, dx_2^{1}
         + dp_2 \,\smwedge \, dq^2 \,\smwedge\, dx_1^{2} \,\smwedge\, dx_2^{2}
         + dp_3 \,\smwedge \, (dq^1+dq^2) \,\smwedge\, dx_1^{3} \,\smwedge\, dx_2^{3} 
\end{equation*}
Any of the three $3$-dimensional subspaces generated by the vectors 
${\partial \over \partial p}$'s or ${\partial \over \partial x_1}$'s or
${\partial \over \partial x_2}$'s is maximal isotropic decomposable with $\mathfrak{N}_L=1$.
To see this, pick the $p$'s subspace, that is, the one generated
by the relation $dx_1=dx_2=dq=0$, and call it $L$. Also, 
denote by $F$ the $3$-isotropic subspace generated 
by $dp=0$. $\omega^\flat (L)$ is spanned by the three $3$-forms on $F$ 
\begin{equation*} 
 \alpha_1 = dq^1 \,\smwedge\, dx_1^{1} \,\smwedge\, dx_2^{1}, \quad 
 \alpha_2 = dq^2 \,\smwedge\, dx_1^{2} \,\smwedge\, dx_2^{2} \quad
 \alpha_3 =(dq^1+dq^2) \,\smwedge\, dx_1^{3} \,\smwedge\, dx_2^{3} \,.
\end{equation*}
\end{exem}
There is no basis $\mathfrak{B}_F^*$ of $F^* \cong L^\bot$ such that, 
if $(\mathfrak{B}_F^*)^n$ is the basis of $\bwedge^n F^*$ generated by 
$\mathfrak{B}_F^*$, 
\footnote{To see this, one can use the fact that the $x_1$'s and $x_2$'s 
subspaces are maximal isotropic decomposable and $3$-dimensional, while their complement 
in $F$ just $2$-dimensional.} 
\[
 \# \, ((\mathfrak{B}_F^*)^n \cap \omega^\flat (L)) =3 \, .
\]  
Therefore $\min_{\mathfrak{B}_F^*} ((\mathfrak{B}_F^*)^n \cap \omega^\flat (L)) = 4$, 
implying $\mathfrak{N}_L =4 - \dim L = 1$.

.

.

%

\textbf{Acknowledgments:} I would like to thank Alan Weinstein and his students
for the opportunity to discuss this subject with and also for their helpful 
comments. Work partially supported by NSF Grant DMS-0707137
               and CNPq (Conselho Nacional de Desenvolvimento
               Cient\'{\i}fico e Tecno\-l\'o\-gico), Brazil


\appendix

\section{Simultaneously integrable Distributions} 
  \label{appendix:Integrability}

 Two distributions $L_1$ and $ L_2$ on a manifold $P$ are simultaneous integrable if 
around each point  $p \in P$ there are integral manifolds $M_1$ of $L_1$ and $M_2$ of $L_2$ 
passing in $p$ such that $M_1 \cap M_2$ is an integral manifold of $L_1 \cap L_2$.

\begin{thm}
 Let $L_1$ and $L_2$ be two integrable distributions such that $\dim (L_1 \cap L_2)$ is 
constant. Then, they are simultaneously 
 integrable if, and only if, $L_1 + L_2$ is integrable.
\end{thm}
\proof
 
We will prove just the ``if'' part and assume that $P$ is a vector space
with $TP=L_1 \oplus L_2$. Let $\{ v_1^{(1)}, \ldots , v_{n_1}^{(1)}, 
u_1, \ldots , u_{n_2}\}$ be a commutative moving frame on $P$
such that $v_1^{(1)}, \ldots , v_{n_1}^{(1)}$ is a basis for $L_1$, which
exists by the integrability hypothesis on $L_1$. Define functions $a_i^j$ 
on $P$ such that  
\[
  e_i^{(2)} = u_i + a_i^j v_j^{(1)} \in L_2 \, .
\]   
Using the fact that $[v_j^{(1)},v_i^{(1)}]=[v_j^{(1)},u_i]=[u_j,u_i]=0$
and $L_2$ is involutive, 
we conclude that $[e_i^{(2)},e_j^{(2)}] \in L_1 \cap L_2=0$, and therefore
\[
  [e_i^{(2)},e_j^{(2)}] = 0 \qquad \text{and} \qquad
  [e_i^{(2)},v] \in L_1 \quad \forall \, v \in L_1 \, .
\]   
In the same way, we can find a basis for $L_1$, and therefore obtain
the moving frame $\{ e_1^{(1)} \ldots , e_{n_1}^{(1)}, 
e_1^{(2)}, \ldots , e_{n_2}^{(2)}\}$ with 
$[e_i^{(1)},e_j^{(2)}] =0$, since it must be in $L_1 \cap L_2=0$.   
\qed


\end{document}